\documentclass[english,12pt]{amsart}
\usepackage{amssymb}
\usepackage{amsfonts}
\usepackage{amscd}

\swapnumbers

\def\Def{{\rm Def}}

\def\11{{\mathbf 1}}
\def\AA{{\mathbf A}}

\def\NN{{\mathbf N}}

\def\PP{{\mathbf P}}
\def\QQ{{\mathbf Q}}
\def\RR{{\mathbf R}}

\def\ZZ{{\mathbf Z}}

\def\cC{{\mathcal C}}

\def\cL{{\mathcal L}}
\def\cM{{\mathcal M}}

\def\cR{{\mathcal R}}
\def\cS{{\mathcal S}}

\def\Sub{{\rm Sub}}

\mathchardef\alphag="7C0B \mathchardef\betag="7C0C
\mathchardef\gammag="7C0D \mathchardef\deltag="7C0E
\mathchardef\varepsilong="7C22 \mathchardef\varphig="7C27
\mathchardef\psig="7C20 \mathchardef\zetag="7C10
\mathchardef\epsilong="7C0F \mathchardef\rhog="7C1A
\mathchardef\taug="7C1C \mathchardef\upsilong="7C1D
\mathchardef\iotag="7C13 \mathchardef\thetag="7C12
\mathchardef\pig="7C19 \mathchardef\sigmag="7C1B
\mathchardef\etag="7C11 \mathchardef\omegag="7C21
\mathchardef\kappag="7C14 \mathchardef\lambdag="7C15
\mathchardef\mug="7C16 \mathchardef\xig="7C18
\mathchardef\chig="7C1F \mathchardef\nug="7C17
\mathchardef\varthetag="7C23 \mathchardef\varpig="7C24
\mathchardef\varrhog="7C25 \mathchardef\varsigmag="7C26
\mathchardef\Omegag="7C0A \mathchardef\Thetag="7C02
\mathchardef\Sigmag="7C06 \mathchardef\Deltag="7C01
\mathchardef\Phig="7C08 \mathchardef\Gammag="7C00
\mathchardef\Psig="7C09 \mathchardef\Lambdag="7C03
\mathchardef\Xig="7C04 \mathchardef\Pig="7C05
\mathchardef\Upsilong="7C07

\newtheorem{prop}[subsubsection]{Proposition}

\theoremstyle{definition}

\newtheorem{def-prop}[subsubsection]{Proposition-Definition}
\newtheorem{def-theorem}[subsubsection]{Theorem-Definition}
\newtheorem{def-lem}[subsubsection]{Lemma-Definition}

\theoremstyle{remark}

\theoremstyle{plain}

\numberwithin{equation}{subsection}

\def\boxit#1#2{\setbox1=\hbox{\kern#1{#2}\kern#1}%
\dimen1=\ht1 \advance\dimen1 by #1 \dimen2=\dp1 \advance\dimen2 by
#1
\setbox1=\hbox{\vrule height\dimen1 depth\dimen2\box1\vrule}%
\setbox1=\vbox{\hrule\box1\hrule}%
\advance\dimen1 by .4pt \ht1=\dimen1 \advance\dimen2 by .4pt
\dp1=\dimen2 \box1\relax}

\def\AA{{\mathbf A}}

\def\NN{{\mathbf N}}

\def\PP{{\mathbf P}}
\def\QQ{{\mathbf Q}}
\def\RR{{\mathbf R}}

\def\ZZ{{\mathbf Z}}

\def\cC{{\mathcal C}}

\def\cL{{\mathcal L}}
\def\cM{{\mathcal M}}

\def\cR{{\mathcal R}}
\def\cS{{\mathcal S}}

\mathchardef\alphag="7C0B \mathchardef\betag="7C0C
\mathchardef\gammag="7C0D \mathchardef\deltag="7C0E
\mathchardef\varepsilong="7C22 \mathchardef\varphig="7C27
\mathchardef\psig="7C20 \mathchardef\zetag="7C10
\mathchardef\epsilong="7C0F \mathchardef\rhog="7C1A
\mathchardef\taug="7C1C \mathchardef\upsilong="7C1D
\mathchardef\iotag="7C13 \mathchardef\thetag="7C12
\mathchardef\pig="7C19 \mathchardef\sigmag="7C1B
\mathchardef\etag="7C11 \mathchardef\omegag="7C21
\mathchardef\kappag="7C14 \mathchardef\lambdag="7C15
\mathchardef\mug="7C16 \mathchardef\xig="7C18
\mathchardef\chig="7C1F \mathchardef\nug="7C17
\mathchardef\varthetag="7C23 \mathchardef\varpig="7C24
\mathchardef\varrhog="7C25 \mathchardef\varsigmag="7C26
\mathchardef\Omegag="7C0A \mathchardef\Thetag="7C02
\mathchardef\Sigmag="7C06 \mathchardef\Deltag="7C01
\mathchardef\Phig="7C08 \mathchardef\Gammag="7C00
\mathchardef\Psig="7C09 \mathchardef\Lambdag="7C03
\mathchardef\Xig="7C04 \mathchardef\Pig="7C05
\mathchardef\Upsilong="7C07

\begin{document}

\title[Real constructible functions]
{Integration of positive constructible functions against Euler
characteristic and dimension}

\author{Raf Cluckers $^*$}\thanks{$^*$ Supported by a postdoctorial fellowship of the Fund for Scientific Research - Flanders (Belgium) (F. W. O) and the European Commission - Marie Curie European Individual Fellowship  HPMF CT 2005-007121.}
\address{Katholieke Universiteit Leuven, Departement wiskunde,
Celestijnenlaan 200B, B-3001 Leu\-ven, Bel\-gium. Current address:
\'Ecole Normale Sup\'erieure, D\'epartement de
ma\-th\'e\-ma\-ti\-ques et applications, 45 rue d'Ulm, 75230 Paris
Cedex 05, France} \email{cluckers@ens.fr}
\urladdr{www.dma.ens.fr/$\sim$cluckers/}

\author{M\'ario Edmundo $^{**}$}\thanks{$^{**}$ With partial support from the FCT 
(Funda\c{c}\~ao para a Ci\^encia e Tecnologia) program POCTI (Portugal/FEDER-EU) 
and Funda\c{c}\~ao Calouste Gulbenkian.}
\address{CMAF Universidade de Lisboa, Avenida Professor Gama Pinto 2, 1649-003 Lisboa, 
Portugal}
\email{edmundo@cii.fc.ul.pt} \urladdr{alf1.cii.fc.ul.pt/$\sim$edmundo}


\begin{abstract}
Following recent work of R.~Cluckers and F.~Loeser [\emph{Fonctions
constructible et int\'egration motivic I}, Comptes rendus de
l'Acad\'emie des Sciences, {\bf 339} (2004) 411 - 416] on motivic
integration, we develop a direct image formalism for positive
constructible functions in the globally subanalytic context. This
formalism is generalized to arbitrary first-order logic models and
is illustrated by several examples on the $p$-adics, on the Presburger structure
and on o-minimal expansions of groups. Furthermore, within this formalism, we define 
the Radon transform and prove the corresponding inversion formula.
\end{abstract}

\maketitle

\renewcommand{\partname}{}

\section{Introduction}
\subsection{}
By a subanalytic set we will always mean a globally subanalytic
subset $X\subset \RR^n$, meaning that $X$ is subanalytic in the classical sense inside 
$\PP^{n}(\RR)$ under the embedding $\RR^n=\AA^{n}(\RR)\subseteq \PP^{n}(\RR)$. 
By a subanalytic function a function
whose graph is a (globally) subanalytic set.

By $\Sub$ we denote the category of subanalytic subsets $X\subset
\RR^n$ for all $n>0$, with subanalytic maps as morphisms. We work with the Euler
characteristic $\chi:\Sub\to\ZZ$ and the dimension
$\dim:\Sub\to\NN$ of subanalytic sets as defined for o-minimal
structures in \cite{vdD}.

Note that if $X\in \Sub$, then  by the o-minimal triangulation theorem in \cite{vdD}, 
the o-minimal Euler characteristic
$\chi (X)$ coincides with the Euler characteristic $\chi _{BM}(X)$ of $X$ with respect 
to the Borel-Moore homology. If
$X\in \Sub$ is locally compact, the o-minimal Euler characteristic $\chi (X)$ coincides 
with the Euler characteristic
$\chi _c(X)$ of $X$ with respect to sheaf cohomology of $X$ with compact supports and 
constant coefficient sheaf.

\subsection{}
By \cite{vdD}, the Euler characteristic $\chi:\Sub\to\ZZ$ satisfies the following
$$\chi (\emptyset)=0,$$
$$\chi (X)=\chi (Y) \,\,{\rm if}\,\,X\,\,{\rm and}\,\,
Y\,\,{\rm are}\,\,{\rm isomorphic}\,\,{\rm in}\,\,\Sub$$
and
$$\chi (X\cup Y)=\chi (X)+\chi (Y)$$
whenever $X,Y\in\Sub$ are disjoint. The last equality for $\chi _{BM}$ and $\chi _c$ 
follows from the long exact
(co)homology sequence. If we take $X$ to be the unit circle in the plane $\RR^2$ and $Y$ 
a point in $X$, we see that
this equality does not hold for the Euler characteristic associated to the topological 
singular (co)homology.

Thus we can think of $\chi:\Sub\to\ZZ$ as a measure with values in the Grothendieck ring 
$K_0(\Sub)$ of the category
$\Sub$ and, for any $X\in \Sub$ and any function $f:X\to \ZZ$ with finite range and the 
property that $f^{-1}(a)\in \Sub$ for
all $a\in \ZZ$ (constructible functions) one has an obvious definition for
$$
\int _{X}f\, \chi
$$
such that $\chi (X)=\int _X1_X\, \chi $  (cf. \cite{viro}).

This measure and integration against Euler characteristic is what is
considered by Viro \cite{viro}, Shapira \cite{schap}, \cite{schap2}
and Brocker \cite{Br}. However, for the measure $\chi:\Sub\to\ZZ$ it
is not true that $\chi (X)=\chi (Y)$  if and only if $X$ and $Y$ are
isomorphic in $\Sub$. Following the recent work of the first author
and Fran\c{c}ois Loeser 
\cite{cr1}, \cite{cr2}, \cite{CL}
on motivic integration, we construct the universal measure $\mu $
for the category $\Sub$ with values in the Grothendieck semi-ring
$SK_0(\Sub)$ of $\Sub$ such that $\mu (X)=\mu (Y)$  if and only if
$X$ and $Y$ are isomorphic in $\Sub$. Furthemore, we develop a
direct image formalism for positive constructible functions, i.e.,
functions $f:X\to SK_0(\Sub)$ with finite range and the property
that $f^{-1}(a)\in \Sub$ for all $a\in SK_0(\Sub)$. This formalism
is generalized to arbitrary first-order logic models and is
illustrated by several examples on the $p$-adics, on the Presburger structure
and on o-minimal expansions of groups. Moreover, within this formalism, we define 
 the Radon transform and prove the corresponding inversion formula.

\section{Positive constructible functions}

We start by pointing out that instead of $\Sub$ we can work in this section with any
o-minimal expansion of a field $R$ using the category $\Def$ whose objects are definable 
sets and whose morphisms are definable maps.

\subsection{}
By a semigroup we mean a commutative monoid with a unit element.
Likewise a semi-ring is a set equipped with two semigroup
structures: addition and multiplication such that $0$ is a unit
element for the addition, $1$ is the unit element for
multiplication and the two operations are connected by $x (y + z)
= xy + xz$ and $0x = 0$. A morphism of semirings is a mapping compatible with the unit
elements and the operations.

\subsection{}
Let $A:=\ZZ\times \NN$ be the semi-ring where addition is given by
$(a,b)+(a',b')=(a+a',\max(b,b'))$, the additive unit element is
$(0,0)$, multiplication is given by $(a,b)(a',b')=(aa',b+b')$ and
the multiplicative unit is $(1,0)$. Note that the ring generated
by $A$ by inverting additively any element of $A$ is $\ZZ$ with
the usual ring structure.

For $Z\in \Sub$ we define $\cC_+(Z)$, as the semi-ring of
functions $Z\to A$ with finite image and whose fibers are
subanalytic sets. We call $\cC_+(Z)$ the semi-ring of positive
constructible functions on $Z$. In particular, $\cC_+(\{0\})=A$.

\subsection{}

If $Z\in
\Sub$ then we denote by $\Sub_Z$ the category of subanalytic maps
$X\to Z$ for $X\in\Sub$ with morphisms subanalytic maps which make
the obvious diagrams commute. We define the Grothendieck semigroup $SK_0 (\Sub_Z)$ as the
quotient of the free abelian semigroup over symbols $[Y
\rightarrow Z]$ with $Y \rightarrow Z$ in $\Sub_Z$ by relations
\begin{equation}\label{eq0}
[\emptyset \rightarrow Z] = 0,
\end{equation}
\begin{equation}\label{eq1}
[Y \rightarrow Z] = [Y' \rightarrow Z]
\end{equation}
if $Y \rightarrow Z$ is isomorphic to $Y' \rightarrow Z$ in
$\Sub_Z$ and
\begin{equation}\label{eq2}
[(Y \cup Y') \rightarrow Z] + [(Y \cap Y') \rightarrow Z] = [Y
\rightarrow Z] + [Y' \rightarrow Z]
\end{equation}
for $Y$ and $Y'$ subsets of some $X \rightarrow Z$. There is a
natural semi-ring structure on $SK_0 (\Sub_Z)$ where the
multiplication is induced by taking fiber products over $Z$.

We write $SK_0 (\Sub)$ for $SK_0 (\Sub_{\{0\}})$ and $[X]$ for $[X\rightarrow \{0\}]$. 
Note that any
element of  $SK_0 (\Sub_Z)$ can be written as $[X\rightarrow Z]$ for some
$X\in\Sub_Z$ because we can take disjoint unions in $\Sub$
corresponding to finite sums in $SK_0 (\Sub_Z)$.

\begin{prop}
For $Z\in\Sub$ there is a natural isomorphism of semi-rings
$$
T:SK_0 (\Sub_Z) \to \cC_+(Z)
$$
induced by sending $[X\to Z]$ in $\Sub_Z$ to $Z\to A:z\mapsto
(\chi(X_z),\dim(X_z))$, where $X_z$ is the fiber above $z$. By
consequence, $SK_0 (\Sub)=A$.
\end{prop}
\begin{proof}
This follows immediately from the trivialisation property for
definable maps in any o-minimal expansion of a field. See \cite{vdD}.
\end{proof}

By means of this result we may identify $SK_0 (\Sub_Z)$ and
$\cC_+(Z)$.

\subsection{Positive measures}
A general notion of positive measures on a Boolean algebra $\cS$
of sets is a map $\mu:\cS\to G$ with $G$ a semigroup satisfying
$$
\mu(X\cup Y)=\mu(X)+\mu(Y)
$$
and
$$\mu(\emptyset)=0$$
whenever $X,Y\in\cS$ are disjoint. Often one has a notion of
isomorphisms between sets in $\cS$ under which the measure should
be invariant and which allows one to take disjoint unions of given
sets in $\cS$ (by taking disjoint isomorphic copies of the sets).

We let $\mu:\Sub\to A$ be the positive measure which sends $X$ to
$(\chi(X),\dim(X))$. This measure is a universal measure on $\Sub$
with the property that $\mu(X)=\mu(Y)$ whenever there exists a
subanalytic bijection between $X$ and $Y$ and where universal
means that any other positive measure with this property
factorises through $\mu$.

Note that $\mu$ measures is in some sense
the topological size since, by the cell decomposition theorem from \cite{vdD}, 
$\mu(A)=\mu(B)$ will hold for two
subanalytic sets $A$, $B$ if and only if for any fixed $n\geq0$
there exists a finite partition of $A$, resp.~$B$ into subanalytic
$C^n$-manifolds $\{A_i\}_{i=1}^m$, resp.~$\{B_i\}_{i=1}^m$ and
subanalytic maps $A_i\to B_i$ which are isomorphisms of
$C^n$-manifolds.

Now we can define the integral of any positive function
$f\in\cC_+(Z)$ as
$$
\int_Z f \, \mu := \sum_i f_i\mu(Z_i)
$$
where $\{Z_i\}$ is any finite partition of $Z$ into subanalytic
sets such that $f$ is constant on each part $Z_i$ with value
$f_i$.

To show that this is independent of the partition $\{Z_i\}$ we just
note that there is a unique $[X\to Z]$ in $SK_0(Z)$ which
corresponds to $f$ under $T$ and that $\sum_i f_i\mu(Z_i)$
corresponds to $[X]=(\chi(X),\dim(X))$ in $A=SK_0(\Sub)$. This
independence follows also from the cell decomposition theorem
(\cite{vdD}).

\subsection{Pushforward}

For $f:X\to Y$ there is an immediate notion of push-forward
$f_!:\cC_+(X)\to \cC_+(Y)$, resp.~$f_!:SK_0(\Sub_X)\to
SK_0(\Sub_Y)$, which is given by
$$
f_!(g)(y)=\int_{f^{-1}(y)} g_{|f^{-1}(y)} \, \mu
$$
for $g\in \cC_+(X)$, resp.~by
$$
f_!([Z\to X])=[Z\to Y],
$$
for $Z\to X$ in $\Sub_X$ and where $Z\to Y$ is given by
composition with $X\to Y$. Note that these pushforwards are
compatible with $T$.

If $Y=\{0\}$, then $SK_0(\Sub_Y)=A$ and we write  $\mu ([Z\to X])$ for
$f_!([Z\to X])$ which is the integral of $[Z\to X]$. Thus the functoriality condition 
$(h\circ f)_{!}=h_{!}\circ f_{!}$ can
be interpreted as Fubini's Theorem, since
$$
\int_{X}g \, \mu=\int_{Y}(\int_{f^{-1}(y)} g_{|f^{-1}(y)} \, \mu )\, \mu
$$
for $g\in \cC_+(X)$ and $h:Y\to \{0\}$.

\subsection{Pullback}
For $f:X\to Y$ a morphism in $\Sub$ there is an immediate notion
of pullback $f^*:\cC_+(Y)\to \cC_+(X)$, resp.~$f^*:SK_0(\Sub_Y)\to
SK_0(\Sub_X)$, which is given by
$$
 f^*(g) = g \circ f
$$
for $g\in \cC_+(Y)$, resp.~by
$$
f^*([Z\to Y])=[Z\otimes_Y X\to X],
$$
for $Z\to Y$ in $\Sub_Y$ and where $Z\otimes_Y X\to X$ is the
projection and $Z\otimes_Y X$ the set-theoretical fiber product.
Note that these pullbacks are also compatible with $T$ and satisfy the functoriality 
property $(f\circ h)^*=h^*\circ f^*$.

\begin{prop}[Projection formula]
Let $f:X\to Y$ be a morphism in $\Sub$ and let $g$ be in
$\cC_+(X)$ and $h$ in $\cC_+(Y)$. Then
$$f_! (g f^* (h)) = f_! (g) h.$$
\end{prop}

\begin{proof}
This is immediate at the level of $SK_0$ since both the
multiplication in $SK_0$ and the pullback are defined by fiber
product.
\end{proof}

\subsection{Radon Transform}
Let $S\subset X\times Y$, $X$, $Y$ be subanalytic sets and write $\pi_X:X\times Y\to X$ 
and  $\pi_Y:X\times Y\to Y$ for the
projections and $q_X=\pi _{X|S}$ and $q_Y=\pi _{Y|S}$. For $g\in \cC_+(X)$, we define the 
Radon transform $\cR _S(g)\in \cC_+(Y)$ by
$$\cR _S(g)=q_{Y!}\circ q_X^*(g)=\pi _{Y!}\circ (\pi _X^*(g)1_S)$$
where $1_S$ is the characteristic function on $S$.

\subsection{Example}
Consider the case $X=\RR^n$, $Y={\rm Gr}(n)$ with $S=\{(p,\Pi):p\in \Pi\}$. Let 
$Z\subseteq \RR^n$ be a subanalytic subset and
$\sigma _Z:{\rm Gr}(n)\to A:\Pi \mapsto (\chi (\Pi\cap Z), {\rm dim}(\Pi \cap Z))$. 
Then $\sigma _Z=\cR_S(1_Z)$.

\medskip
Let $S'\subset Y\times X$ be another subanalytic set and put
$q'_X=\pi _{X|S'}$ and $q'_Y=\pi _{Y|S'}$. The following proposition
is proved just like in \cite{schap2}.

\begin{prop}[Inversion formula]\label{prop inv}
Let $r:S\otimes_Y S' \to X\times X$ be the projection and suppose that the following
 hypothesis hold

\vspace{.05in}
\noindent (*) there exists  $\lambda \in A$ such that
$[r^{-1}(x,x')]=\lambda $ for all $x\neq x'$, $x,x'\in X$;

\vspace{.05in}
\noindent (**) there exists $0\neq \theta \in A$ such that 
$[r^{-1}(x,x)]=\theta +\lambda $ for all $x\in X$.

\vspace{.05in}
\noindent If $g$ is in $\cC_+(X)$, then
\begin{equation}\label{eq inv radon}
\cR_{S'}\circ \cR_{S}(g) = \theta g+ \lambda \int _X g\, \mu 
\end{equation}
and this is independent of the choice of $\theta $.
\end{prop}

\begin{proof}
Let $h$ and $h'$ be the projections from $S\otimes_Y S'$ to $S$ and $S'$ respectively.
Then by definition of fiber product, $q_Y\circ h=q'_Y\circ h'$, and so 
by functoriality of pullback and 
pushforward we have $h'_!\circ h^*=q'^*_Y\circ q_{Y!}$. Thus
$\cR_{S'}\circ \cR _S(g)=q'_{X!}\circ (q'_{Y})^*\circ q_{Y!}\circ q_X^*(g)
=q'_{X!}\circ h'_!\circ h^*\circ q_{X}^*(g)$.

The last formula is also equal to $p_{2!}\circ r_{!}\circ r^*\circ p_{1}^*(g)$ where 
$p_1,p_2:X\times X\to X$ are the
projections onto the first and second coordinates respectively, since 
$q_X\circ h=p_1\circ r$ and $q'_X\circ h'=p_2\circ r$. The hypothesis shows 
that
$r_!(1_{S\otimes_Y S'})=\theta 1_{\Delta _X}+\lambda 1_{X\times X}$, moreover, this 
expression is independent of the choice
of $\theta $. By the projection formula,
$r_!(r^*(p^*_1(g)))=r_!(1_{S\otimes _Y S'}r^*(p^*_1(g)))=r_!(1_{S\otimes _Y S'})p^*_1(g)
=(\theta 1_{\Delta _X}+\lambda 1_{X\times X})p^*_1(g)$ holds, 
hence we obtain
$p_{2!}((\theta 1_{\Delta _X}+\lambda 1_{X\times X})p^*_1(g))
=\theta p_{2!}(1_{\Delta _X} p^*_1(g))+\lambda p_{2!}(p^*_1(g))
=\theta g+ \lambda \int _X g\, \mu $ as required.

We now show that the inversion formula is independent of the choice of $\theta$. If 
$\theta +\lambda =\theta '+\lambda $ and
$\theta \neq \theta '$, then necessarily $\lambda _2>\theta _2$, $\lambda _2>\theta '_2$ 
and 
$\theta _1=\theta '_1$ with
$\lambda =(\lambda _1,\lambda _2)$, $\theta =(\theta _1,\theta _2)$ and 
$\theta '=(\theta '_1,\theta '_2)$. Hence,
$\theta g+\lambda \int _X g\, \mu =\theta 'g+\lambda \int _X g\, \mu $ for all 
$x\in X$.
\end{proof}

\subsection{Example}
Consider the case $X=\RR^n$, $Y={\rm Gr}(n)$ with $S=\{(p,\Pi):p\in \Pi\}$ and 
$S'=\{(\Pi,p):p\in \Pi\}$. Then
$[r^{-1}(x,x)]=[\PP^{n-1}]$ and $[r^{-1}(x,x')]=[\PP^{n-2}]$ for all $x,x'\in \RR^n$ 
with $x\neq x'$. Since
$[\PP ^n]=(\frac{1+(-1)^n}{2},n)$, we have
$$\cR_{S'}\circ \cR_{S}(g) = ((-1)^{n+1},n-1)g+ (\frac{1+(-1)^n}{2},n-2)\int _Xg\, \mu .$$

In particular, we have
$$\cR_{S'}\circ \cR_{S}(1_Z) = ((-1)^{n+1},n-1)1_Z+ (\frac{1+(-1)^n}{2},n-2)[Z]$$
for every subanalytic subset $Z$ of $\RR^n$ .

\section{Direct image formalism in model theory}
Let $\cM$ be a model of a theory in a language $\cL$ with at least two
constant symbols $c_1,c_2$ satisfying $c_1\not=c_2$. For $Z$ a
definable set we define the category $\Def_Z(\cM)$, also written
$\Def_Z$ for short, whose objects are definable sets $X$ with a
definable map $X\to Z$ and whose morphisms are definable
maps which make the obvious diagram commute. We write
$\Def(\cM)$ or $\Def$ for $\Def_{\{c_1\}}(\cM)$. In $\cM$ one can pursue the
usual operations of set theory like finite unions, intersections,
Cartesian products, disjoint unions and fiber products.

We define the Grothendieck semigroup $SK_0 (\Def_Z)$ as the
quotient of the free abelian semigroup over symbols $[Y
\rightarrow Z]$ with $Y \rightarrow Z$ in $\Def_Z$ by relations
\begin{equation}
[\emptyset \rightarrow Z] = 0,
\end{equation}
\begin{equation}
[Y \rightarrow Z] = [Y' \rightarrow Z]
\end{equation}
if $Y \rightarrow Z$ is isomorphic to $Y' \rightarrow Z$ in
$\Def_Z$ and
\begin{equation}
[(Y \cup Y') \rightarrow Z] + [(Y \cap Y') \rightarrow Z] = [Y
\rightarrow Z] + [Y' \rightarrow Z]
\end{equation}
for $Y$ and $Y'$ subsets of some $X \rightarrow Z$. There is a
natural semi-ring structure on $SK_0 (\Def_Z)$ where the
multiplication is induced by taking fiber products over $Z$. Note
that any element of  $SK_0 (\Def_Z)$ can be written as $[X\to Z]$
for some $X\to Z\in\Def_Z$ because we can take disjoint unions in $\cM$
corresponding to finite sums in $SK_0 (\Def_Z)$.

The map $\Def\to SK_0 (\Def)$ sending $X$ to its class $[X]$ is a
universal positive measure with the property that two sets have
the same measure if there exists a definable bijection between
them. For $f:X\to Y$ there is an immediate notion of push-forward
$f_!:SK_0(\Def _X)\to SK_0(\Def _Y)$ given by
$$
f_!([Z\to X])=[Z\to Y],
$$
for $Z\to X$ in $\Def _X$ and where $Z\to Y$ is given by
composition with $X\to Y$.

If $Y=\{c_1\}$, then we write  $\mu ([Z\to X])$ for
$f_!([Z\to X])$ which  we call the integral of $[Z\to X]$, note that  $\mu ([Z\to X])$ is just
 $[Z]$ in $SK_0(\Def)$. Thus the functoriality 
condition $(f\circ h)_{!}=f_{!}\circ h_{!}$
can be interpreted as Fubini's Theorem.

There is also an immediate notion of
pullback $f^*:SK_0(\Def _Y)\to SK_0(\Def _X)$ given by
$$
f^*([Z\to Y])=[Z\otimes_Y X\to X],
$$
for $Z\to Y$ in $\Def _Y$ and where $Z\otimes_Y X\to X$ is the
projection and $Z\otimes_Y X$ the set-theoretical fiber product. The pullback is 
functorial, i.e., $(f\circ h)^*=h^*\circ f^*$.

\begin{prop}[Projection formula]
Let $f:X\to Y$ be a morphism in $\Def$ and let $g$ be in
$SK_0(\Def _X)$ and $h$ in $SK_0(\Def _Y)$. Then
$$f_! (g f^* (h)) = f_! (g) h.$$
\end{prop}
\begin{proof}
Exactly the same proof as for the subanalytic sets above
works.
\end{proof}

\subsection{Radon transform and inversion formula}
One can also define the Radon transform in this context in exactly the same way as in the subanalytic 
case. Furthermore, the same argument as in the subanalytic case gives the corresponding inversion 
formula. However, since in general there is no trivialisation theorem, the conditions (*) and (**) in 
Proposition (\ref{prop inv}) have to be replaced by global conditions.Using the embedding 
$SK_0(\Def)\to SK_0(\Def _U)$ sending $[W]$ to $[W\times U\to U]$ where $W\times U\to U$ is the 
projection, 
 the statement becames: 

\vspace{.05in}
\noindent Let $r:S\otimes_Y S' \to X\times X$ be the projection and suppose that the following
 hypothesis hold

\vspace{.05in}
\noindent (*) there exists $Z_1$ in $\Def$ such that in $SK_0(\Def _{X_1})$ we have 
$$[B_1\to X_1]=[Z_1],$$

\vspace{.05in}
\noindent (**) there exists $Z_2$ in $\Def$ such that in $SK_0(\Def _{\Delta _X})$ we have 
$$[B_2\to \Delta _X]=[Z_1]+[Z_2]$$ 

\vspace{.05in}
\noindent 
where $X_1=X\times X\setminus \Delta _X$, $B_1=S\otimes _YS'\setminus r^{-1}(\Delta _X)$, 
$B_2=S\otimes _YS'\cap r^{-1}(\Delta _X)$ and $B_1\to X_1$ and $B_2\to \Delta _X$ are the restrictions 
of the projection $r:S\otimes _YS'\to X\times X'$.
If $Z\to X$ is in $\Def_X$, then
\begin{equation}\label{eq inv mod radon}
\cR_{S'}\circ \cR_{S}([Z\to X]) = [Z_2][Z\to X]+ [Z_1][Z] 
\end{equation}
and this is independent of the choice of $Z_2$. 


\section{Examples}

\subsection{Semialgebraic and subanalytic sets in  $\QQ_p$}
For $K$ any finite field extension of the field $\QQ_p$ of $p$-adic
numbers, one can calculate explicitely the semi-ring of
semialgebraic sets  $SK_0(K,\mathrm{Sem})$, resp.~of globally
subanalytic sets $SK_0(K,\Sub)$, using work of \cite{C} for
semialgebraic sets, resp.~using work of \cite{Ccell} for the
subanalytic sets.
In both cases it is a subset of $\NN\times\NN$ and the class of a
semialgebraic set $X$, resp.~a subanalytic set $X$, is $(\sharp
X,0)$ if $X$ is finite and $(0,\dim X)$ if $X$ is infinite. This is
because there exists a semialgebraic bijection between two infinite
semialgebraic sets if and only if they have the same dimension, and
similarly for subanalytic sets. However, no trivialisation theorem
is known hence the relative semi-Grothendieck rings
$SK_0(K,\mathrm{Sem}_Z)$, resp.~$SK_0(K,\Sub_Z)$ for $Z$
semialgebraic, resp.~subanalytic, are expected to be much more
complicated than maps $Z\to \NN\times\NN$ with finite image.

\subsection{Presburger sets}
Consider the Presburger structure on $\ZZ$ by using the Presburger
language
$$
\cL_{\rm PR} = \{+, -, 0, 1, \leq\} \cup \{\equiv_n\ \mid n\in \NN,\
n
>1\},
$$
with $\equiv_n$ the equivalence relation modulo $n$. Again one can
can calculate explicitely the semi-ring $SK_0(\ZZ,\cL_{\rm PR})$,
using work of \cite{pres}. It is a subset of $\NN\times\NN$ and the
class of a Presburger set $X$ is $(\sharp X,0)$ if $X$ is finite and
$(0,\dim X)$ if $X$ is infinite, where the dimension of \cite{pres}
is used. Again this is because there exists a Presburger bijection
between two infinite Presburger sets if and only if they have the
same dimension. Again, no trivialisation theorem is known hence the
relative semi-Grothendieck rings are expected to be more
complicated.

\subsection{Semilinear sets}\label{expl sl}
Let $K=(K,0,1,+,\cdot ,<)$ be an ordered field and consider the structure 
$\cM =(K,0,1,+,(\lambda _c)_{c\in K},<)$, where
$\lambda _c$ is the scalar multiplication by $c\in K$. The category $\Def$ in this case 
is the category of $K$-semilinear
sets with $K$-semilinear maps.

By \cite{Mari}, the Grothendieck ring $K_0(\Def)$ is isomorphic to $E=\ZZ[x]/(x(x+1))$ 
and there is a universal Euler
characteristic $\epsilon :\Def \to E$ (see also \cite{KF}).

Let $D$ be the set whose elements are of the form 
$\sum _{i=1}^ny^{k_i}z^{l_i}\in \NN[y,z]$ with $k_i\leq l_i$,  and for
$i\neq j$, $\neg (y^{k_i}z^{l_i}=y^{k_j}z^{l_j})\wedge 
\neg (y^{k_i}z^{l_i}\prec y^{k_j}z^{l_j})\wedge
\neg (y^{k_j}z^{l_j}\prec y^{k_i}z^{l_i})$. Here, $y^{k_i}z^{l_i}\prec y^{k_j}z^{l_j}$ 
if and only if $k_i<k_j$ and $l_i<l_j$.

The set $D$ can be equipped with a semi-ring structure in the following way: the zero 
element $0_D$ is
$\sum _{i=1}^0y^{k_i}z^{l_i}$, the identity element $1_D$ is $y^0z^0$, the addition is 
given by
$$\sum _{i=1}^ny^{k_i}z^{l_i}+_D\sum _{i=1}^my^{k'_i}z^{l'_i}
=\sum {\rm max}_{\prec}\{y^kz^l:y^kz^l \,{\rm a}\,
{\rm monomial}\,{\rm in}\,\sum _{i=1}^ny^{k_i}z^{l_i}+\sum _{i=1}^my^{k'_i}z^{l'_i}\}
$$
and multiplication by
$$\sum _{i=1}^ny^{k_i}z^{l_i}\cdot _D\sum _{i=1}^my^{k'_i}z^{l'_i}
=\sum {\rm max}_{\prec}\{y^kz^l:y^kz^l \,{\rm a}\,
{\rm monomial}\,{\rm in}\,\sum _{i=1}^ny^{k_i}z^{l_i}\cdot \sum _{i=1}^my^{k'_i}z^{l'_i}\}
$$
where the symbol $\sum {\rm max}_{\prec}S$ mean that we sum up the $\prec $-maximal 
elements of the finite set $S$.

By \cite{Mari}, there is a universal abstract dimension $\delta :\Def\to D$ and two 
sets in $\Def$ are isomorphic in
$\Def$ if and only if they have the same universal Euler characteristic and the same 
universal abstract dimension. Thus,
if $A$ is the semi-ring $E\times D$, then the Grothendieck semi-ring $SK_0(\Def)$ is 
isomorphic to $A$ and the map
$\mu :\Def\to A$ given by $\mu (X)=(\epsilon (X), \delta (X))$ is the positive universal 
measure on $\Def$.

Note that the results we used above from \cite{Mari} were proved in the field of real 
numbers but the same arguments
hold in any arbitrary ordered field $K$.

\subsection{Semi-bounded sets}
Let $K=(K,0,1,+,\cdot ,<)$ be a real closed field and consider the structure 
$\cM =(K,0,1,+,(\lambda _c)_{c\in K},B,<)$,
where $\lambda _c$ is the scalar multiplication by $c\in K$ and $B$ is the graph of 
multiplication on a bounded interval.
The category $\Def$ in this case is the category of $K$-semibounded sets with 
$K$-semibounded maps.

By \cite{MPP} all bounded semialgebraic subsets are in $\Def$ and, by \cite{PSS}, 
$\cM$ is, up to definability, the only
o-minimal structure properly between $(K,0,1,+,(\lambda _c)_{c\in K},<)$ and 
$(K,0,1,+,\cdot ,<)$.

By  \cite{Mari}, the Grothendieck ring $K_0(\Def)$ is isomorphic to $E=\ZZ[x]/(x(x+1))$ 
and there is a universal Euler
characteristic $\epsilon :\Def \to E$ (see also \cite{KF}). Furthermore, if $D$ is the semi-ring of 
Example \ref{expl sl}, then there is a universal
abstract dimension $\delta :\Def\to D$ and two sets in $\Def$ are isomorphic in 
$\Def$ if and only if they have the same
universal Euler characteristic and the same universal abstract dimension. Thus, if 
$A$ is the semi-ring $E\times D$, then
the Grothendieck semi-ring $SK_0(\Def)$ is isomorphic to $A$ and the map 
$\mu :\Def\to A$ given by
$\mu (X)=(\epsilon (X), \delta (X))$ is the positive universal measure on $\Def$.

The results we used above from \cite{Mari} were proved in the field of real numbers 
and are based on Peterzil's \cite{Pe}
structure theorem for semibounded sets in the real numbers. However, the same arguments 
hold in any arbitrary real closed
field $K$ using the structure theorem from \cite{Ed}.

\bibliographystyle{amsplain}

\begin{thebibliography}{SGA}



%
%
%
%
\bibitem{Br}
L. Br\"{o}cker, \textit{Euler integration and Euler multiplication}, Adv. Geom., \textbf{5} (2005), 145-169.
%
%
%

\bibitem{C} R. Cluckers,
\textit{Classification of semialgebraic $p$-adic sets up to semialgebraic bijection}, J. 
Reine Angew. Math.,
\textbf{540} (2001), 105-114.

\bibitem{pres}
R. Cluckers, \textit{Presburger sets and p-minimal fields}, J.
Symbolic Logic, \textbf{68} (2003), 153--162.

\bibitem{Ccell} R. Cluckers, \textit{Analytic $p$-adic cell decomposition and integrals}, 
Trans. Amer. Math. Soc.,
\textbf{356} (2004), 1489 - 1499. Available
at arXiv:math.NT/0206161.

\bibitem{cr1}{R. Cluckers,  F. Loeser},
\textit{Fonctions constructibles et int\'egration motivique I}, C.
R. Math. Acad. Sci. Paris  \textbf{339}  (2004),   411--416.


\bibitem{cr2}{R. Cluckers,  F. Loeser},
\textit{Fonctions constructibles et int\'egration motivique II}, C.
R. Math. Acad. Sci. Paris  \textbf{339}  (2004), 487--492.


\bibitem{CL}
R. Cluckers, F. Loeser, \textit{Constructible motivic functions and
motivic integrals}, Available
at arXiv:math.AG/0410203.

%
%
%
%
%
%
%
%
%
%
%
%
%
%
%
%
%
%
%
%
%
%
%
%
\bibitem{vdD}
L. van den Dries, \textit{Tame topology and o-minimal structures},
Cambridge University Press, Lecture note series, \textbf{248},
(1998).


%
%
%
\bibitem{Ed}
M. Edmundo, \textit{Structure theorems for o-minimal expansions of groups}, Ann. Pure 
Appl. Logic, \textbf{102} (2000), 159--181.

%
%
%
%
%
\bibitem{KF}
M. Kageyama, M. Fujita, \textit{Grothendieck rings of o-minimal expansions 
of ordered abelian groups}, to appear in J. Algebra. Available
at arXiv:math.LO/0505341.

%
\bibitem{MPP}
D. Marker, Y. Peterzil, A. Pillay, \textit{Additive reducts of real closed fields}, J. 
Symbolic Logic, \textbf{57} (1992), 109--117.
%
\bibitem{Mari}
J. Ma\u{r}ikov\'a, \textit{Geometric properties of definable sets}, Diplomov\'a Pr\'ace. 
Univerzita Karlova v Praze, Praha, 2003.
%
%
\bibitem{Pe}
Y. Peterzil \textit{A structure theorem for semibounded sets in the reals}, J. Symbolic 
Logic \textbf{57} (1992), 779--794.
%
\bibitem{PSS}
A. Pillay, P. Scowcroft, C.Steinhorn, \textit{Between groups and rings}, Rocky Mountain J. 
Math. \textbf{9} (1989), 871--885.
%

\bibitem{schap}
P. Schapira, \textit{Operations on constructible functions}, J.
Pure Appl. Algebra, \textbf{72} (1991), 83--93.

\bibitem{schap2}
P. Schapira, \textit{Tomography of constructible functions} In: AAECC, G.Cohen, M.Giusti 
and T.Mora (eds.) Lecture Notes in
Computer Sciences, Springer Verlag, \textbf{948} (1995), 427--435.
%
%


\bibitem{viro}
O. Viro, \textit{Some integral calculus based on Euler
characteristic}, Lecture Notes in Math, Springer-Verlag,
\textbf{1346} (1988), 127--138.



\end{thebibliography}

\end{document}